\documentclass[12pt,reqno]{amsart} 

\usepackage{mathrsfs,amsfonts,amsthm,amsmath,amssymb,thmtools} 
\usepackage[dvipsnames]{xcolor} 

\usepackage{enumitem}                
\usepackage{lineno}                 

\usepackage[numbers,sort&compress]{natbib}  

\usepackage{graphics}
\usepackage{graphicx}   
\graphicspath{{images/}} 
\usepackage{subfig}     
\usepackage{rotating}

\usepackage[hyperindex, pdfencoding = auto, psdextra, bookmarksdepth = 4]{hyperref}  
\hypersetup{
    bookmarksopen      = true,
    bookmarksopenlevel = 1,
    bookmarksnumbered  = true,
    pdftitle    = {},
    pdfauthor   = {},
    linktoc     = page,
    anchorcolor = green,
    colorlinks  = true,        
    linkcolor   = black,
    citecolor   = blue,
    filecolor   = blue,
    urlcolor    = magenta,
}

\usepackage[capitalise, noabbrev, nameinlink]{cleveref}    

\theoremstyle{plain}

\newtheorem{theorem}{Theorem}[section]
\newtheorem{lemma}[theorem]{Lemma}

\newtheorem{proposition}[theorem]{Proposition}
\newtheorem{corollary}[theorem]{Corollary}
\newtheorem{question}[theorem]{Question}

\theoremstyle{definition}   
\newtheorem{definition}[theorem]{Definition}

\makeatletter
    \def\subsection{\@startsection{subsection}{2}%
    \z@{.5\linespacing\@plus.7\linespacing}{.3\linespacing}%
    {\normalfont\bfseries}}
\makeatother 

\makeatletter                                         
    \newcommand{\LeftEqNo}{\let\veqno\@@leqno}        
\makeatother

\numberwithin{equation}{section}                      

\begin{document}

\

\vspace{-2cm}

\title[Cleanness and real rank zero in $C^*$-algebras]
{Cleanness and real rank zero in $C^*$-algebras}

\author{Lu Cui}
\address{Lu Cui, College of Science, China University of Petroleum (East China), Qingdao, 266580, China}
\email{nilu@upc.edu.cn}

\author{Minghui Ma}
\address{Minghui Ma, School of Mathematical Sciences, Dalian University of Technology, Dalian, 116024, China}
\email{minghuima@dlut.edu.cn}

\author{Wei Yuan}
\address{Wei Yuan, Institute of Mathematics, Academy of Mathematics and Systems Science, Chinese Academy
of Sciences, Beijing, 100190, China; School of Mathematical Sciences, University of Chinese
Academy of Sciences, Beijing 100049, China}
\email{wyuan@math.ac.cn}


\makeatletter
\@namedef{subjclassname@2020}{\textup{2020} Mathematics Subject Classification}
\makeatother

\begin{abstract}
A unital $C^*$-algebra is called clean if every element is the sum of an idempotent and an invertible element.
We prove that cleanness is equivalent to real rank zero.
More precisely, the invertible summand can be chosen with inverse norm at most $2+\varepsilon$ for every $\varepsilon>0$, and the constant $2$ is optimal.
In particular, every von Neumann algebra is clean, answering a question of Lam.
Clean decompositions with uniformly bounded idempotent summands characterize real rank zero and stable rank one.
Moreover, $*$-cleanness, defined by requiring the idempotent summands to be projections, is equivalent to real rank zero, stable rank one, and triviality of the $K_1$-group.
Finally, sharp cleanness, defined by requiring clean decompositions with inverse norm at most $2$, is equivalent to real rank zero and the $SAW^*$ property.
\end{abstract}

\maketitle

\section{Introduction}

A unital $C^*$-algebra has \emph{real rank zero} \cite{BP-91} if its invertible self-adjoint elements are norm dense in its self-adjoint part.
Recall that a unital ring $R$ is said to be an \emph{exchange ring} \cite[Lemma~1.1]{AGOP-98} if for every element $a\in R$, there exists an idempotent $e\in R$ such that $e\in Ra$ and $1-e\in R(1-a)$.
Ara, Goodearl, O'Meara, and Pardo \cite[Theorem~7.2]{AGOP-98} proved that a unital $C^*$-algebra is an exchange ring if and only if it has real rank zero.

A unital ring is called \emph{clean} if every element can be expressed as the sum of an idempotent and an invertible element.
Nicholson \cite[Proposition~1.8(1)]{Nicholson-77} introduced this notion in his study of exchange rings and proved that every clean ring is an exchange ring.
Consequently, every clean unital $C^*$-algebra has real rank zero.
\v{S}ter \cite[Corollary~3.7]{Ster-14} proved that every unital $C^*$-algebra of real rank zero is a corner of a clean complex algebra and left open whether such $C^*$-algebras are themselves clean.
Our first result answers this question affirmatively and provides clean decompositions with a uniform bound on the inverse norms of the invertible summands.

\begin{theorem}\label{thm clean}
Let $\mathcal A$ be a unital $C^*$-algebra.
The following statements are equivalent:
\begin{enumerate}
\item[$(1)$] $\mathcal A$ has real rank zero;
\item[$(2)$] $\mathcal A$ is clean;
\item[$(3)$] for every $T\in\mathcal A$ and every $\varepsilon>0$, there exists an idempotent $P\in\mathcal A$ such that $T-P$ is invertible and $\|(T-P)^{-1}\|\leqslant2+\varepsilon$.
\end{enumerate}
\end{theorem}

The main new implication in the preceding theorem is $(1)\Rightarrow(3)$.
Together with \cite[Theorem~7.2]{AGOP-98} and \cite[Example~4.3(1)]{Ster-14}, \Cref{thm clean} shows that cleanness, the exchange property, and weak cleanness in the sense of \v{S}ter are equivalent for unital $C^*$-algebras.
Since real rank zero passes to hereditary $C^*$-subalgebras, \Cref{thm clean} implies that cleanness passes to projection corners of unital $C^*$-algebras.
For general rings, corners of clean rings need not be clean \cite{Ster-12}.

In 2005, Lam asked which von Neumann algebras are clean as rings; see \cite[Introduction]{Vas-10}.
Akalan and Va\v s \cite[Corollary~2.10]{AV-13} proved that every finite $AW^*$-algebra is almost clean, meaning that every element is the sum of an idempotent and an element that is neither a left nor a right zero-divisor.
It was proved in \cite{CHWYZ-22} that all finite von Neumann algebras and all infinite factors acting on a separable Hilbert space are clean.
These results were subsequently extended in \cite{CM-26} to all finite $AW^*$-algebras and all countably decomposable infinite $AW^*$-factors.
Since every $AW^*$-algebra has real rank zero, \Cref{thm clean} shows that all $AW^*$-algebras, and hence all von Neumann algebras, are clean.
This answers Lam's question, as well as \cite[Question~5.12]{CM-26} and \cite[Question~6.3]{CHWYZ-22}.
\Cref{thm clean} also implies that $\mathcal A\otimes_{\min}\ell^\infty(\mathbb N)$ is clean whenever $\mathcal A$ is a clean unital $C^*$-algebra.
This gives a negative answer to \cite[Question~6.4]{CHWYZ-22}.

We next consider bounds on the norms of the idempotent summands in clean decompositions.
Recall that a unital $C^*$-algebra has \emph{stable rank one} \cite{Rieffel-83} if its invertible elements are norm dense.
Our second result asserts that a unital $C^*$-algebra has real rank zero and stable rank one if and only if it admits clean decompositions with uniformly bounded idempotents.

\begin{theorem}\label{thm bounded-idempotents}
Let $\mathcal A$ be a unital $C^*$-algebra.
The following statements are equivalent:
\begin{enumerate}
\item[$(1)$] $\mathcal A$ has real rank zero and stable rank one;
\item[$(2)$] there exists $\varepsilon>0$ such that for every $T\in\mathcal A$, there exists an idempotent $P\in\mathcal A$ such that $T-P$ is invertible and $\|P\|\leqslant1+\varepsilon$;
\item[$(3)$] for every $0<\varepsilon\leqslant1$ and every $T\in\mathcal A$, there exists an idempotent $P\in\mathcal A$ such that $T-P$ is invertible, $\|P\|\leqslant1+\varepsilon$, and $\|(T-P)^{-1}\|\leqslant17\varepsilon^{-1/2}$.
\end{enumerate}
\end{theorem}

Since an idempotent in a $C^*$-algebra is a projection if and only if its norm is at most $1$, it is natural to ask when the bound $1$ can be attained.
Although \Cref{thm bounded-idempotents} allows the idempotent norms to approach $1$, the corresponding inverse norm bound tends to infinity as $\varepsilon\to0$.
The existence of projection summands therefore requires
a separate argument.
Recall that a unital $C^*$-algebra is called \emph{$*$-clean} if every element is the sum of a projection and an invertible element \cite{Vas-10}.
We denote by $K_1(\mathcal A)$ the first topological $K$-group of $\mathcal A$; see \cite{Blackadar-98}.
The following theorem characterizes $*$-cleanness in terms of real rank, stable rank, and $K$-theory.

\begin{theorem}\label{thm star-clean}
Let $\mathcal A$ be a unital $C^*$-algebra.
The following statements are equivalent:
\begin{enumerate}
\item[$(1)$] $\mathcal A$ has real rank zero, stable rank one, and $K_1(\mathcal A)=\{0\}$;
\item[$(2)$] $\mathcal A$ is $*$-clean;
\item[$(3)$] for every $T\in\mathcal A$, there exists a projection $P\in\mathcal A$ such that $T-P$ is invertible and $\|(T-P)^{-1}\|\leqslant10^5$;
\item[$(4)$] there exists $\varepsilon>0$ such that for every $T\in\mathcal A$, there exists an idempotent $P\in\mathcal A$ such that $T-P$ is invertible, $\|P\|\leqslant1+\varepsilon$, and $\|(T-P)^{-1}\|\leqslant\frac12\varepsilon^{-1/2}$.
\end{enumerate}
\end{theorem}

As an application, we show in \Cref{cor AW-star-clean} that an $AW^*$-algebra is $*$-clean if and only if it is finite.
This answers \cite[Question~6.1]{CHWYZ-22} and \cite[Question~4.4]{CM-26}.

The estimates in \Cref{thm bounded-idempotents,thm star-clean} reveal a quantitative role for $K_1(\mathcal A)$.
For a unital $C^*$-algebra $\mathcal A$ of real rank zero and stable rank one, the order $\varepsilon^{-1/2}$ is optimal as $\varepsilon\to0$ when $K_1(\mathcal A)\ne\{0\}$.
When $K_1(\mathcal A)=\{0\}$, the idempotent summands can be chosen to be projections, with a universal bound on the inverse norms of the invertible summands.
For related connections between exchange rings, real rank zero, and $K$-theory, see Ara, Goodearl, O'Meara, and Raphael \cite{AGOR-00} and Perera \cite{Perera-00}.

We finally return to the inverse norm estimate in \Cref{thm clean}.
The constant $2$ in condition $(3)$ of \Cref{thm clean} is optimal.
Indeed, for $T=(1/2)I$ and any idempotent $P$, the identity $(I-2P)^2=I$ shows that
\begin{equation*}
    \|(T-P)^{-1}\|=2\|I-2P\|\geqslant2.
\end{equation*}
It is therefore natural to ask whether $\varepsilon$ can be omitted from condition $(3)$.
This motivates the following definition.

\begin{definition}\label{def sharp}
A unital $C^*$-algebra $\mathcal A$ is \emph{sharply clean} if for every $T\in\mathcal A$, there exists an idempotent $P\in\mathcal A$ such that $T-P$ is invertible and $\|(T-P)^{-1}\|\leqslant2$.
\end{definition}

Every $AW^*$-algebra is sharply clean by \Cref{prop AW}.
In particular, the optimal inverse norm bound for von Neumann algebras is $2$, answering \cite[Question~6.2]{CHWYZ-22}.

Recall that a $C^*$-algebra $\mathcal A$ is called an \emph{$SAW^*$-algebra} if for any positive $A,B\in\mathcal A$ with $AB=0$, there exists a positive $X\in\mathcal A$ such that $XA=A$ and $XB=0$.
Pedersen \cite{Pedersen-86} introduced this notion as a noncommutative analogue of sub-Stonean spaces.

\begin{theorem}\label{thm sharply-clean}
Let $\mathcal A$ be a unital $C^*$-algebra.
The following statements are equivalent:
\begin{enumerate}
\item[$(1)$] $\mathcal A$ has real rank zero and is an $SAW^*$-algebra;
\item[$(2)$] $\mathcal A$ is sharply clean.
\end{enumerate}
\end{theorem}

As a consequence, we show in \Cref{cor sharply-clean-consequences} $(1)$ that a $\sigma$-finite unital $C^*$-algebra is sharply clean if and only if it is an $AW^*$-algebra.
In particular, every sharply clean unital $C^*$-algebra admitting a faithful state is an $AW^*$-algebra.
Since every separable $SAW^*$-algebra is finite dimensional by \cite[Corollary~2]{Pedersen-86}, every separable sharply clean unital $C^*$-algebra is finite dimensional; see \Cref{cor sharply-clean-consequences} $(2)$.
Consequently, every infinite-dimensional separable unital $C^*$-algebra of real rank zero is clean but is not sharply clean.
In particular, the term $\varepsilon$ in \Cref{thm clean} cannot be omitted in general.
By \Cref{thm star-clean} and \Cref{cor sharply-clean-consequences} $(2)$, every separable
infinite-dimensional unital AF algebra is $*$-clean but is not
sharply clean.

The characterizations in \Cref{thm clean,thm bounded-idempotents,thm star-clean} show that cleanness, $*$-cleanness, and the existence of clean decompositions with uniformly bounded idempotent summands are preserved under matrix amplification.
Han and Nicholson \cite{HN-01} proved that matrix rings over clean unital rings are clean.
It is therefore natural to ask whether sharp cleanness has the same permanence property.
Pedersen \cite[p.~175]{Pedersen-89} considered whether the $SAW^*$ property is preserved under passage to matrix algebras.
The abelian case follows from Smith and Williams \cite[Theorem~3.4]{SW-88}.
Since real rank zero is preserved under matrix amplification, \Cref{thm sharply-clean} identifies the following question with the real rank zero case of Pedersen's problem.

\begin{question}\label{qu matrix-sharply-clean}
Let $\mathcal A$ be a sharply clean unital $C^*$-algebra.
Is $M_n(\mathcal A)$ sharply clean for every integer $n\geqslant2$?
\end{question}

The paper is organized as follows.
In \Cref{sec clean}, we prove \Cref{thm clean,thm bounded-idempotents} and show that every $AW^*$-algebra is sharply clean.
In \Cref{sec star-clean}, we prove \Cref{thm star-clean} and characterize $*$-clean $AW^*$-algebras.
In \Cref{sec sharply-clean}, we prove \Cref{thm sharply-clean} and derive consequences for $\sigma$-finite and separable $C^*$-algebras, as well as corona algebras.

\section{Cleanness and real rank zero}
\label{sec clean}

We begin with a standard inverse norm estimate and include an elementary proof for completeness.

\begin{lemma}\label{lem norm-inverse}
Let $\mathcal A$ be a unital $C^*$-algebra.
If $T\in\mathcal A$ satisfies $\operatorname{Re}T\geqslant cI$ for some $c>0$, then $T$ is invertible with $\|T^{-1}\|\leqslant c^{-1}$.
\end{lemma}

\begin{proof}
Since $\operatorname{Re}T\geqslant cI$, we have
\begin{equation*}
\begin{aligned}
    T^*T
    &=(T-cI)^*(T-cI)+2c\operatorname{Re}T-c^2I
      \geqslant c^2I,\\
    TT^*
    &=(T-cI)(T-cI)^*+2c\operatorname{Re}T-c^2I
      \geqslant c^2I.
\end{aligned}
\end{equation*}
Thus, both $T^*T$ and $TT^*$ are invertible in $\mathcal A$.
Consequently, $T$ is also invertible.
Moreover, we have
\begin{equation*}
    \|T^{-1}\|^2
    =\|T^{-1}(T^{-1})^*\|
    =\|(T^*T)^{-1}\|
    \leqslant c^{-2},
\end{equation*}
which gives $\|T^{-1}\|\leqslant c^{-1}$.
This completes the proof.
\end{proof}

The following lemma constructs a clean decomposition with an explicit inverse norm bound from diagonal estimates for $\operatorname{Re}T$.

\begin{lemma}\label{lem idempotent-construction}
Let $\mathcal A$ be a unital $C^*$-algebra.
Suppose that $T\in\mathcal A$ and there are projections $E,F\in\mathcal A$ with $I=E+F$ such that
\begin{equation*}
	E(\operatorname{Re}T)E\leqslant(1/2+\delta)E,\qquad F(\operatorname{Re}T)F\geqslant(1/2-\delta)F
\end{equation*}
for some $0\leqslant\delta<1/2$.
Then there exists an idempotent $P$ in $\mathcal A$ such that $T-P$ is invertible and
\begin{equation*}
	\|P\|\leqslant 1+2\|T\|,\qquad\|(T-P)^{-1}\|\leqslant\frac{1}{1/2-\delta}=\frac{2}{1-2\delta}.
\end{equation*}
\end{lemma}

\begin{proof}
Let $X=E(T-T^*)F$ and $P=E+X$.
Since $EF=FE=0$, we have
\begin{equation*}
    EX=X,\qquad XE=0,\qquad X^2=0.
\end{equation*}
Consequently, $P^2=P$.
Let $J=F-E$.
Then $J=J^*=J^{-1}$.
With respect to the decomposition $I=E+F$, we have
\begin{equation*}
    S:=J(T-P)
    =
    \begin{pmatrix}
        E-ETE & -ET^*F\\
        FTE   & FTF
    \end{pmatrix}.
\end{equation*}
It follows that
\begin{equation*}
    \operatorname{Re}S
    =
    \begin{pmatrix}
        E-E(\operatorname{Re}T)E & 0\\
        0 & F(\operatorname{Re}T)F
    \end{pmatrix}
    \geqslant(1/2-\delta)I.
\end{equation*}
By \Cref{lem norm-inverse}, $S$ is invertible and
$\|S^{-1}\|\leqslant(1/2-\delta)^{-1}$. Since $J$ is a self-adjoint unitary operator, $T-P=JS$ is invertible and
\begin{equation*}
    \|(T-P)^{-1}\|
    =\|S^{-1}\|
    \leqslant\frac{1}{1/2-\delta}
    =\frac{2}{1-2\delta}.
\end{equation*}
It is clear that $\|P\|\leqslant 1+2\|T\|$.
This completes the proof.
\end{proof}

Cleanness of finite $AW^*$-algebras and countably decomposable infinite $AW^*$-factors was established in \cite{CM-26}.
The following proposition extends these results to all $AW^*$-algebras and gives the optimal inverse norm bound $2$, improving the previous bound $4$ in \cite[Theorem~4.2]{CM-26}.

\begin{proposition}\label{prop AW}
Let $\mathcal A$ be an $AW^*$-algebra.
Then for every element $T\in\mathcal A$, there exists an idempotent $P\in\mathcal A$ such that $T-P$ is invertible and
\begin{equation*}
	\|(T-P)^{-1}\|\leqslant2.
\end{equation*}
In particular, $\mathcal A$ is sharply clean.
\end{proposition}

\begin{proof}
Let $H=\operatorname{Re}T$.
By \cite[Lemma~2.5]{CM-26}, there are projections $E,F\in\mathcal A$ such that
\begin{equation*}
    E+F=I,\qquad
    EHE\leqslant\frac12E,\qquad
    FHF\geqslant\frac12F.
\end{equation*}
The conclusion follows from \Cref{lem idempotent-construction} with $\delta=0$.
\end{proof}

In a $C^*$-algebra of real rank zero, approximation by invertible self-adjoint elements gives the required inequalities in \Cref{lem idempotent-construction} with arbitrarily small $\delta>0$.
This yields the following approximate version of \Cref{prop AW}.

\begin{lemma}\label{lem RR0-clean}
Let $\mathcal A$ be a unital $C^*$-algebra of real rank zero.
Then for every $T\in\mathcal A$ and every $\varepsilon>0$, there exists an idempotent $P$ in $\mathcal A$ such that $T-P$ is invertible and
\begin{equation*}
	\|P\|\leqslant 1+2\|T\|,\qquad\|(T-P)^{-1}\|\leqslant2+\varepsilon.
\end{equation*}
In particular, $\mathcal A$ is clean.
\end{lemma}

\begin{proof}
Let $H=\operatorname{Re}T$ and $\delta=\varepsilon/(2(2+\varepsilon))$.
Then $0<\delta<1/2$.
Since $\mathcal A$ has real rank zero, there exists an invertible self-adjoint element $S$ in $\mathcal A$ such that
\begin{equation*}
    \|S-(H-(1/2)I)\|<\delta.
\end{equation*}
Since $0\notin\sigma(S)$, continuous functional calculus gives the projections
\begin{equation*}
    E=\chi_{(-\infty,0)}(S)\in\mathcal A,\qquad F=I-E.
\end{equation*}
Then $ESE\leqslant0$ and $FSF\geqslant0$.
It follows that
\begin{equation*}
\begin{aligned}
    EHE
    &\leqslant ESE+(1/2+\delta)E
     \leqslant(1/2+\delta)E,\\
    FHF
    &\geqslant FSF+(1/2-\delta)F
     \geqslant(1/2-\delta)F.
\end{aligned}
\end{equation*}
By \Cref{lem idempotent-construction}, there exists an
idempotent $P$ in $\mathcal A$ such that $T-P$ is invertible and
\begin{equation*}
    \|(T-P)^{-1}\|
    \leqslant\frac{2}{1-2\delta}
    =2+\varepsilon.
\end{equation*}
Since $T=P+(T-P)$, it follows that $\mathcal A$ is clean.
\end{proof}

We next prepare the quantitative construction used in the proof of \Cref{thm bounded-idempotents}.

\begin{lemma}\label{lem finite-spectrum-bounded}
Let $\mathcal A$ be a unital $C^*$-algebra of stable rank one, let $U\in\mathcal A$ be unitary, and let $H\in\mathcal A$ be positive with finite spectrum.
Then for every $0<t\leqslant1$, there exists an idempotent $P\in\mathcal A$ such that $UH-P$ is invertible and
\begin{equation*}
    \|P\|^2\leqslant1+t^2,\qquad
    \|(UH-P)^{-1}\|\leqslant\frac{16}{t}.
\end{equation*}
\end{lemma}

\begin{proof}
Let
\begin{equation*}
	T=UH,\qquad E=\chi_{[0,r)}(H),\qquad F=I-E,\qquad 0<r<1.
\end{equation*}
If $E=0$, take $P=0$; then $T$ is invertible with $\|T^{-1}\|\leqslant r^{-1}$.
If $F=0$, take $P=I$; then $T-I$ is invertible with $\|(T-I)^{-1}\|\leqslant(1-r)^{-1}$.
Suppose that $0<E<I$.
With respect to the decomposition $I=E+F$, we write
\begin{equation*}
	U=
	\begin{pmatrix}
		U_{11} & U_{12}\\
		U_{21} & U_{22}
	\end{pmatrix},\qquad
	H=
	\begin{pmatrix}
		H_1 & 0\\
		0 & H_2
	\end{pmatrix}.
\end{equation*}
Then $\|H_1\|\leqslant r$ and $\|H_2^{-1}\|\leqslant r^{-1}$.
Let
\begin{equation*}
	B=U_{11}^*-H_1,\qquad A=(BB^*)^{1/2}.
\end{equation*}
Since $U_{11}^*U_{11}+U_{21}^*U_{21}=E$ and $H_1^2\leqslant r^2E$, we have
\begin{align*}
	BB^*+U_{21}^*U_{21}
	&=E-rU_{11}^*U_{11}-\frac{1-r}{r}H_1^2+r(U_{11}-r^{-1}H_1)^*(U_{11}-r^{-1}H_1)\\
	&\geqslant(1-r)^2E.
\end{align*}
Moreover, $\|A\|\leqslant1+r$, and hence $BB^*=A^2\leqslant(1+r)A$.
It follows that
\begin{equation*}
    A+tU_{21}^*U_{21}
    \geqslant\frac{t}{1+r}(BB^*+U_{21}^*U_{21})
    \geqslant\alpha tE,\qquad \alpha=\frac{(1-r)^2}{1+r}.
\end{equation*}
By \cite[Proposition~4.1]{Blackadar-04}, $E\mathcal AE$ has stable rank one.
Thus, there are invertible elements $\{B_n\}_{n=1}^\infty$ in $E\mathcal AE$ converging to $B$.
Write their polar decompositions as
\begin{equation*}
	B_n=(B_nB_n^*)^{1/2}W_n,
\end{equation*}
where each $W_n$ is unitary in $E\mathcal AE$.
Then $\|B-AW_n\|\to0$.
Let $0<\beta<\alpha$.
Then $\|B-AW\|<\beta t$ for some unitary $W\in E\mathcal AE$.
Let $P=E+tU_{21}W$.
It is straightforward to verify that $P^2=P$, $\|P\|^2\leqslant1+t^2$, and
\begin{equation*}
    U^*(T-P)=
    \begin{pmatrix}
        K&0\\
        L&H_2
    \end{pmatrix},
    \qquad
    \begin{cases}
        K=-B-tU_{21}^*U_{21}W,\\
        L=-U_{12}^*-tU_{22}^*U_{21}W.
    \end{cases}
\end{equation*}
Since $-KW^*=A+tU_{21}^*U_{21}+(B-AW)W^*$, we obtain
\begin{equation*}
	\operatorname{Re}(-KW^*)\geqslant(\alpha-\beta)tE.
\end{equation*}
By \Cref{lem norm-inverse}, $K$ is invertible in $E\mathcal AE$ with $\|K^{-1}\|\leqslant\bigl((\alpha-\beta)t\bigr)^{-1}$.
We also have
\begin{equation*}
	L=-FU^*PE,\qquad \|L\|\leqslant\|P\|\leqslant\sqrt{1+t^2}\leqslant\sqrt2.
\end{equation*}
Therefore, $U^*(T-P)$ is invertible with inverse
\begin{equation*}
	\begin{pmatrix}
		K^{-1} & 0\\
		-H_2^{-1}LK^{-1} & H_2^{-1}
	\end{pmatrix},
\end{equation*}
and hence
\begin{align*}
	\|(T-P)^{-1}\|
    &\leqslant\max\{\|K^{-1}\|,\|H_2^{-1}\|\}+\|H_2^{-1}LK^{-1}\|\\
    &\leqslant\frac1t\left(\max\left\{\frac1{\alpha-\beta},\frac1r\right\}+\frac{\sqrt2}{r(\alpha-\beta)}\right).
\end{align*}
Let $r=1/3$.
Then $\alpha=1/3$.
Since $3+9\sqrt2<16$, we may choose $\beta>0$ sufficiently small such that
\begin{equation*}
    \|(T-P)^{-1}\|
    \leqslant\frac{1+3\sqrt2}{(1/3-\beta)t}
    \leqslant\frac{16}{t}.
\end{equation*}
The same bound also holds in the cases $E=0$ and $F=0$ considered above.
This completes the proof.
\end{proof}

We can now prove \Cref{thm clean,thm bounded-idempotents}.

\begin{proof}[Proof of \Cref{thm clean}]
The implication $(3)\Rightarrow(2)$ is immediate, and $(1)\Rightarrow(3)$ follows from \Cref{lem RR0-clean}.
Finally, $(2)\Rightarrow(1)$ follows from \cite[Proposition~1.8(1)]{Nicholson-77} and \cite[Theorem~7.2]{AGOP-98}.
\end{proof}

\begin{proof}[Proof of \Cref{thm bounded-idempotents}]
The implication $(3)\Rightarrow(2)$ is immediate.

$(2)\Rightarrow(1)$.
Let $\varepsilon>0$ be as in $(2)$.
Since $\mathcal A$ is clean, it has real rank zero by \Cref{thm clean}.
For every $T\in\mathcal A$ and every $n\geqslant1$, there exists an idempotent $P_n\in\mathcal A$ such that $nT-P_n$ is invertible and $\|P_n\|\leqslant1+\varepsilon$.
Then $T-n^{-1}P_n$ is invertible and converges in norm to $T$.
Thus, $\mathcal A$ has stable rank one.

$(1)\Rightarrow(3)$.
Let $T\in\mathcal A$, $0<\varepsilon\leqslant1$, and $t=\sqrt{\varepsilon}$.
By stable rank one and real rank zero, elements of the form $UH$, with $U$ unitary and $H$ positive of finite spectrum, are norm dense in $\mathcal A$.
Choose such an element $UH$ with $\|T-UH\|<t/272$, and let $P$ be the idempotent provided by \Cref{lem finite-spectrum-bounded}.
Then $\|P\|\leqslant\sqrt{1+\varepsilon}\leqslant1+\varepsilon$ and
\begin{equation*}
    \|(T-P)^{-1}\|
    \leqslant\frac{16/t}{1-(16/t)\|T-UH\|}
    <17\varepsilon^{-1/2}.
\end{equation*}
This completes the proof.
\end{proof}

\section{\texorpdfstring{$*$-cleanness and the $K_1$-group}
{*-cleanness and the K1-group}}
\label{sec star-clean}

We first recall two facts about unitary elements.
Write $\mathcal U(\mathcal A)$ for the unitary group of $\mathcal A$ and $\mathcal U_0(\mathcal A)$ for the connected component of $I$.
By \cite[Theorem~2.10]{Rieffel-87}, if $\mathcal A$ has
stable rank one, then the canonical map
\begin{equation*}
    \mathcal U(\mathcal A)/\mathcal U_0(\mathcal A)
    \to K_1(\mathcal A)
\end{equation*}
is an isomorphism.
Moreover, if $\mathcal A$ has real rank zero, then every unitary element in $\mathcal U_0(\mathcal A)$ is a norm limit of unitary elements with finite spectrum by \cite{Lin-93}.

The following application of \Cref{lem norm-inverse} provides the inverse norm estimate needed for the proof of \Cref{thm star-clean}.

\begin{lemma}\label{lem unitary-positive}
Let $\mathcal A$ be a unital $C^*$-algebra.
Suppose that $H\in\mathcal A$ is a positive element and $U\in\mathcal A$ is a unitary element with finite spectrum.
Then there exists a self-adjoint unitary $J\in\mathcal A$ such that $UH+J$ is invertible and
\begin{equation*}
    \|(UH+J)^{-1}\|\leqslant1+\|H\|.
\end{equation*}
\end{lemma}

\begin{proof}
Write $U=\sum_{j=1}^n\lambda_jE_j$, where $\{E_j\}_{j=1}^n$ are mutually orthogonal projections with sum $I$ and $\lambda_j\in\mathbb T$.
Let $J=\sum_{j=1}^n\varepsilon_jE_j$, where $\varepsilon_j\in\{-1,1\}$ such that $\varepsilon_j\operatorname{Re}\lambda_j\geqslant0$ for each $1\leqslant j\leqslant n$.
Then
\begin{equation*}
    J=J^*=J^{-1},\qquad
    \operatorname{Re}(JU)
    =\sum_{j=1}^n\varepsilon_j\operatorname{Re}\lambda_jE_j
    \geqslant0.
\end{equation*}
Let $K=H^{1/2}$ and $T=I+KJUK$.
Then $\operatorname{Re}T=I+K\operatorname{Re}(JU)K\geqslant I$.
It follows from \Cref{lem norm-inverse} that $T^{-1}\in\mathcal A$ and $\|T^{-1}\|\leqslant1$.
Recall that
\begin{equation*}
	(I+BA)^{-1}=I-B(I+AB)^{-1}A
\end{equation*}
whenever $I+AB$ is invertible.
By taking $A=K$ and $B=JUK$, we have
\begin{equation*}
	(I+JUH)^{-1}=I-JUKT^{-1}K.
\end{equation*}
Since $UH+J=J(I+JUH)$ and $J$ is unitary, $UH+J$ is invertible and
\begin{equation*}
    \|(UH+J)^{-1}\|
    =\|(I+JUH)^{-1}\|
    \leqslant1+\|K\|^2
    =1+\|H\|.
\end{equation*}
This completes the proof.
\end{proof}

Let $\mathcal A$ be a unital $C^*$-algebra of real rank zero with $K_1(\mathcal A)=\{0\}$.
It is well known that $K_1(E\mathcal AE)=\{0\}$ for every nonzero projection $E\in\mathcal A$.
Indeed, let $\mathcal I$ be the closed ideal of $\mathcal A$ generated by $E$.
It follows from \cite[Theorems~2.10 and~3.14]{BP-91} that $K_0(\mathcal A)\to K_0(\mathcal A/\mathcal I)$ is surjective.
The six-term exact sequence therefore shows that $K_1(\mathcal I)\to K_1(\mathcal A)$ is injective, and hence $K_1(\mathcal I)=\{0\}$.
Since $E$ is full in $\mathcal I$, we have $K_1(E\mathcal AE)\cong K_1(\mathcal I)=\{0\}$ by Morita invariance of $K$-theory.

We next derive the following consequence of \Cref{lem unitary-positive}.

\begin{lemma}\label{lem corner-symmetry}
Let $\mathcal A$ be a unital $C^*$-algebra of real rank zero and stable rank one with $K_1(\mathcal A)=\{0\}$.
Then for every projection $E\in\mathcal A$, every $r>0$, and every $Y\in E\mathcal AE$ with $\|Y\|<r$, there exists a self-adjoint unitary $J$ in $E\mathcal AE$ such that $Y+J$ is invertible in $E\mathcal AE$ and
\begin{equation*}
    \|(Y+J)^{-1}\|<1+r.
\end{equation*}
\end{lemma}

\begin{proof}
The assertion is immediate if $E=0$.
If $E\ne 0$, then $E\mathcal AE$ has real rank zero and stable rank one by \cite{BP-91,Blackadar-04} and $K_1(E\mathcal AE)=\{0\}$.
By \cite{Lin-93} and \cite[Theorem~2.10]{Rieffel-87}, every unitary element in $E\mathcal AE$ is a norm limit of unitary elements with finite spectrum.
Together with stable rank one and polar decomposition, this gives a sequence $\{U_nH_n\}_{n=1}^\infty$ converging to $Y$, where each $U_n\in E\mathcal AE$ is unitary with finite spectrum and each $H_n\in E\mathcal AE$ is positive.
By \Cref{lem unitary-positive}, there exists a self-adjoint unitary $J_n\in E\mathcal AE$ such that $U_nH_n+J_n$ is invertible and
\begin{equation*}
    \|(U_nH_n+J_n)^{-1}\|\leqslant1+\|H_n\|.
\end{equation*}
Since $\|H_n\|=\|U_nH_n\|\to\|Y\|<r$, we have
\begin{equation*}
    \|(Y+J_n)^{-1}\|
    \leqslant
    \frac{1+\|H_n\|}
    {1-(1+\|H_n\|)\|Y-U_nH_n\|}
    <1+r
\end{equation*}
for all sufficiently large $n$.
This completes the proof.
\end{proof}

We can now prove \Cref{thm star-clean}.

\begin{proof}[Proof of \Cref{thm star-clean}]
The implication $(3)\Rightarrow(2)$ is immediate, and $(3)\Rightarrow(4)$ follows by taking $\varepsilon=(2\cdot10^5)^{-2}$.

$(4)\Rightarrow(2)$.
Let $\varepsilon>0$ be as in $(4)$.
It is clear that $\varepsilon\leqslant1/16$ by applying $(4)$ to the element $(1/2)I$.
Let $T\in\mathcal A$, let $P\in\mathcal A$ be an idempotent as in $(4)$, and let $Q\in\mathcal A$ be the range projection of $P$.
Then $QP=P$ and $PQ=Q$.
If $P=0$, then $T$ is invertible and the conclusion holds.
Suppose that $P\ne0$.
With respect to the decomposition $I=Q+(I-Q)$, we write
\begin{equation*}
    P=\begin{pmatrix}Q&X\\0&0\end{pmatrix}.
\end{equation*}
Then $\|P-Q\|=\|X\|$ and $\|P\|^2=\|PP^*\|=\|Q+XX^*\|=1+\|X\|^2$.
It follows that
\begin{equation*}
    \|P-Q\|^2=\|P\|^2-1
    \leqslant(1+\varepsilon)^2-1
    =2\varepsilon+\varepsilon^2,
\end{equation*}
and hence
\begin{equation*}
    \|(T-P)^{-1}\|\,\|P-Q\|
    \leqslant\frac12\sqrt{2+\varepsilon}<1.
\end{equation*}
Thus, $T-Q$ is invertible.
Therefore, $\mathcal A$ is $*$-clean.

$(2)\Rightarrow(1)$.
Suppose that $\mathcal A$ is $*$-clean.
For every $T\in\mathcal A$ and every $\varepsilon>0$, applying $*$-cleanness to $\varepsilon^{-1}T$ gives a projection $P\in\mathcal A$ such that $T-\varepsilon P$ is invertible.
Since $\|\varepsilon P\|\leqslant\varepsilon$ and $T-\varepsilon P$ is self-adjoint whenever $T$ is, $\mathcal A$ has stable rank one and real rank zero.
For every $U\in\mathcal U(\mathcal A)$, applying $*$-cleanness to $(I+U)/2$ gives a projection $P\in\mathcal A$ such that $U+J$ is invertible, where $J=I-2P$.
It follows that $-1\notin\sigma(JU)$ and hence $JU\in\mathcal U_0(\mathcal A)$ by continuous functional calculus.
Since $J=e^{i\pi P}\in\mathcal U_0(\mathcal A)$, we have
\begin{equation*}
	U=J(JU)\in\mathcal U_0(\mathcal A).
\end{equation*}
Thus, $\mathcal U(\mathcal A)=\mathcal U_0(\mathcal A)$, and therefore $K_1(\mathcal A)=\{0\}$ by \cite[Theorem~2.10]{Rieffel-87}.

$(1)\Rightarrow(3)$.
Let $F(t)=t(2t+1)^2$ for $t>0$.
Suppose that $A\in\mathcal A$ is invertible, $E\in\mathcal A$ is a nonzero projection, $\lambda\in\{i,-i\}$, and $t>\|A^{-1}\|$.
We claim that there exists a self-adjoint unitary $J\in E\mathcal AE$ such that
\begin{equation}\label{equ A+J}
    \|(A+J-\lambda E)^{-1}\|<F(t).
\end{equation}
Let $B=A^{-1}$.
Then $\|B\|<t$.
By \Cref{lem corner-symmetry}, applied to $Y=\lambda E+2EBE$ with $\|Y\|<1+2t$, there exists a self-adjoint unitary $J\in E\mathcal AE$ such that
\begin{equation*}
    K=(\lambda E+2EBE+J)^{-1},
    \qquad \|K\|<2+2t.
\end{equation*}
Since $(J-\lambda E)^{-1}=(J+\lambda E)/2$ in $E\mathcal AE$, it is routine to verify that
\begin{equation*}
    (A+J-\lambda E)^{-1}=B-2BEKEB.
\end{equation*}
This proves the claim \eqref{equ A+J}.
Now let $T\in\mathcal A$, $X=2T-I$, and $t>\sqrt2$.
As in the proof of \Cref{lem corner-symmetry}, we may choose a positive $H\in\mathcal A$ and a unitary $U\in\mathcal A$ with finite spectrum such that
\begin{equation*}
    \|X-HU\|<\frac1{\sqrt2}-\frac1t.
\end{equation*}
Write $U=\sum_{j=1}^nu_jE_j$.
For every $1\leqslant j\leqslant n$, we can choose $v_j\in\{1,-1,i,-i\}$ such that $\operatorname{Re}(v_j\overline{u_j})\geqslant1/\sqrt2$.
Let $V=\sum_{j=1}^nv_jE_j$.
Then $V$ is unitary and
\begin{equation*}
    \sigma(V)\subseteq\{1,-1,i,-i\},
    \qquad
    \operatorname{Re}(VU^*)
    =\sum_{j=1}^n
    \operatorname{Re}(v_j\overline{u_j})E_j
    \geqslant\frac1{\sqrt2}I.
\end{equation*}
It follows that
\begin{equation*}
    \operatorname{Re}((X+V)U^*)
    \geqslant H+
    \left(\frac1{\sqrt2}-\|X-HU\|\right)I
    >t^{-1}I.
\end{equation*}
By \Cref{lem norm-inverse}, $X+V$ is invertible and $\|(X+V)^{-1}\|<t$.

Let $E_+=\chi_{\{i\}}(V)$ and $E_-=\chi_{\{-i\}}(V)$.
Applying \eqref{equ A+J} with $A=X+V$, $E=E_+$, and $\lambda=i$, we obtain a self-adjoint unitary $J_+\in E_+\mathcal AE_+$ such that
\begin{equation*}
    \|(X+V+J_+-iE_+)^{-1}\|<F(t).
\end{equation*}
If $E_+=0$, then we take $J_+=0$.
Applying \eqref{equ A+J} again, with $A=X+V+J_+-iE_+$, $E=E_-$, $\lambda=-i$, and $F(t)$ in place of $t$, we obtain a self-adjoint unitary $J_-\in E_-\mathcal AE_-$ such that
\begin{equation*}
    \|(X+V+J_+-iE_++J_-+iE_-)^{-1}\|<F(F(t)).
\end{equation*}
If $E_-=0$, then we take $J_-=0$.
Since $E_+$ and $E_-$ are orthogonal spectral projections of $V$, the element
\begin{equation*}
    S=(V-iE_++iE_-)+(J_++J_-)
\end{equation*}
satisfies $S=S^*$ and $S^2=I$.
Thus, $P=(I-S)/2$ is a projection and
\begin{equation*}
    \|(T-P)^{-1}\|
    =2\|(X+S)^{-1}\|
    <2F(F(t)).
\end{equation*}
Since $2F(F(\sqrt2))=37\,024+26\,658\sqrt2<10^5$, we complete the proof by taking $t>\sqrt2$ sufficiently close to $\sqrt2$.
\end{proof}

By \cite[Theorem~4.3]{CM-26}, every $*$-clean $AW^*$-algebra is finite.
The following corollary establishes the equivalence of $*$-cleanness and finiteness.

\begin{corollary}\label{cor AW-star-clean}
An $AW^*$-algebra is $*$-clean if and only if it is finite.
In particular, a von Neumann algebra is $*$-clean if and only if it is finite.
\end{corollary}

\begin{proof}
Let $\mathcal A$ be an $AW^*$-algebra.
It is well known that $\mathcal A$ has real rank zero and $K_1(\mathcal A)=\{0\}$.
If $\mathcal A$ is $*$-clean, then it has stable rank one by \Cref{thm star-clean} and is therefore finite.
Conversely, suppose that $\mathcal A$ is finite.
Then every $T\in\mathcal A$ admits a unitary polar decomposition $T=U|T|$ by \cite[Lemma~2]{Handelman-79}.
Since the invertible elements $U(|T|+\varepsilon I)$ converge in norm to $T$ as $\varepsilon\to0$, $\mathcal A$ has stable rank one.
Thus, the conclusion follows from \Cref{thm star-clean}.
\end{proof}

\section{Sharp cleanness and \texorpdfstring{$SAW^*$-algebras}{SAW*-algebras}}
\label{sec sharply-clean}

We first record an auxiliary consequence of sharp cleanness that will be used in the proof of \Cref{thm sharply-clean}.

\begin{lemma}\label{lem sharp-anticommutator}
Let $\mathcal A$ be a sharply clean unital $C^*$-algebra.
For every self-adjoint contraction $H\in\mathcal A$, there exists a self-adjoint contraction $K\in\mathcal A$ such that
\begin{equation*}
    HK+KH\geqslant I-K^2\geqslant0.
\end{equation*}
\end{lemma}

\begin{proof}
By sharp cleanness applied to $(I+H)/2$, there exists an idempotent $P\in\mathcal A$ such that
\begin{equation*}
    X=(I+H-2P)^{-1},\qquad \|X\|\leqslant1.
\end{equation*}
Since $(X^{-1}-H)^2=(I-2P)^2=I$, we obtain
\begin{equation*}
    XH+HX=I-X^2+XH^2X.
\end{equation*}
Let $K=\operatorname{Re}X$, $B=\operatorname{Im}X$, and $D=I-H^2$.
Then $0\leqslant D\leqslant I$ and $\|K\|\leqslant1$.
Taking real parts gives
\begin{equation*}
    HK+KH
    =I-KDK+BDB
    \geqslant I-K^2
    \geqslant0.
\end{equation*}
This completes the proof.
\end{proof}

We identify $\mathcal A$ with its canonical image in its double dual $\mathcal A^{**}$, which is the enveloping von Neumann algebra of $\mathcal A$.
We next derive the following approximate separation result from \Cref{lem sharp-anticommutator}.

\begin{lemma}\label{lem approximate-separation}
Let $\mathcal A$ be a unital sharply clean $C^*$-algebra, let $\{E_n\}_{n=1}^\infty$ and $\{F_n\}_{n=1}^\infty$ be projections in $\mathcal A$ such that all these projections are mutually orthogonal, and let
\begin{equation*}
    E=\sum_{n=1}^\infty E_n,\qquad
    F=\sum_{n=1}^\infty F_n,\qquad G=E+F,
\end{equation*}
where the sums converge strongly in $\mathcal A^{**}$.
Then for every $\varepsilon>0$, there exists a self-adjoint element $K\in\mathcal A$ such that
\begin{equation*}
    \|K\|\leqslant1,\qquad KG=GK,\qquad\|GKG-(E-F)\|<\varepsilon.
\end{equation*}
\end{lemma}

\begin{proof}
Let $Q_{2n-1}=F_n$, $Q_{2n}=E_n$, and $S=E-F$.
For $0<t<1$, let
\begin{equation*}
    H=\sum_{n=1}^\infty(-t^2)^nQ_n.
\end{equation*}
The series converges in norm, $\|H\|\leqslant t^2$, and the kernel projection of $H$ is $I-G$.
By \Cref{lem sharp-anticommutator}, there exists a self-adjoint contraction $K\in\mathcal A$ such that
\begin{equation*}
    0\leqslant I-K^2\leqslant HK+KH\leqslant2t^2I.
\end{equation*}
Let $D=HK+KH$.
Since $(I-G)D(I-G)=0$ and $D\geqslant0$, we have $D(I-G)=0$.
Thus, $HK(I-G)=0$, which gives $GK(I-G)=0$ and hence $KG=GK$.
For any $n\ne m$, positivity of $D$ gives
\begin{align*}
	\bigl|(-t^2)^n+(-t^2)^m\bigr|\|Q_nKQ_m\|
	&=\|Q_nDQ_m\|\\
	&\leqslant\|Q_nDQ_n\|^{1/2}\|Q_mDQ_m\|^{1/2}
	\leqslant2t^{n+m}.
\end{align*}
It follows that
\begin{equation*}
    \|Q_nKQ_m\|
    \leqslant\frac{2t^{|n-m|}}{1-t^{2|n-m|}}\qquad(n\ne m).
\end{equation*}
Let $L=\sum_{n=1}^\infty Q_nKQ_n$, where the sum converges strongly in $\mathcal A^{**}$.
Then $\|L\|\leqslant1$.
Each row and column of the preceding bounds has sum at most
\begin{equation*}
    4\sum_{k=1}^\infty\frac{t^k}{1-t^{2k}}
    \leqslant\frac{4t}{(1-t)^2}.
\end{equation*}
Let $\delta_t=\frac{4t}{(1-t)^2}$.
The Schur estimate gives $\|KG-L\|\leqslant\delta_t$.
Moreover, $Q_nDQ_n\geqslant0$ implies that $(-1)^nQ_nKQ_n\geqslant0$.
Since $SQ_n=Q_nS=(-1)^nQ_n$ and $\|K\|\leqslant1$, we obtain
\begin{equation*}
    0\leqslant SL=LS
    =\sum_{n=1}^\infty(-1)^nQ_nKQ_n
    \leqslant G.
\end{equation*}
Since $S^2=G$ and $GL=L$, we have $(SL)^2=L^2$ and hence $0\leqslant G-SL\leqslant G-L^2$.
It follows that
\begin{equation*}
	\|S-L\|=\|G-SL\|\leqslant\|G-L^2\|\leqslant\|G-K^2G\|+\|(KG)^2-L^2\|\leqslant2t^2+2\delta_t.
\end{equation*}
Therefore, $\|KG-(E-F)\|\leqslant2t^2+3\delta_t\to0$ as $t\to0$.
Taking $t>0$ sufficiently small completes the proof.
\end{proof}

We can now prove \Cref{thm sharply-clean}.

\begin{proof}[Proof of \Cref{thm sharply-clean}]
$(1)\Rightarrow(2)$.
Let $T\in\mathcal A$ and $H=\operatorname{Re}T$.
By \cite[Theorem~2.2]{DZ-03}, there exists a projection $F\in\mathcal A$ such that
\begin{equation*}
    \chi_{(1/2,\infty)}(H)\leqslant F\leqslant\chi_{[1/2,\infty)}(H).
\end{equation*}
Let $E=I-F$.
It is clear that $EHE\leqslant(1/2)E$ and $FHF\geqslant(1/2)F$.
By \Cref{lem idempotent-construction} with $\delta=0$, there exists an idempotent $P\in\mathcal A$ such that $T-P$ is invertible and $\|(T-P)^{-1}\|\leqslant2$.
Thus, $\mathcal A$ is sharply clean.

$(2)\Rightarrow(1)$.
Since $\mathcal A$ is clean, it has real rank zero by \Cref{thm clean}.
Let $A,B\in\mathcal A$ be nonzero positive elements with $AB=0$, and let $E,F\in\mathcal A^{**}$ be their support projections.
Then $EF=0$.
The hereditary subalgebras generated by $A$ and $B$ are $\sigma$-unital and have real rank zero.
By \cite[Corollary~2.8 and Proposition~2.9]{BP-91}, they admit increasing approximate identities of projections.
Taking successive differences gives mutually orthogonal projections $\{E_n\}_{n=1}^\infty$ and $\{F_n\}_{n=1}^\infty$ in $\mathcal A$ such that
\begin{equation*}
    \sum_{n=1}^\infty E_n=E,\qquad
    \sum_{n=1}^\infty F_n=F
\end{equation*}
strongly in $\mathcal A^{**}$.
Let $G=E+F$ and $\mathcal B=\{X\in\mathcal A\colon XG=GX\}$.
The map
\begin{equation*}
    \pi\colon\mathcal B\to G\mathcal A^{**}G,
    \qquad \pi(X)=GXG,
\end{equation*}
is a unital $*$-homomorphism and has closed range.
Hence $E-F\in\pi(\mathcal B)$ by \Cref{lem approximate-separation}.
Since $EF=0$, it follows that $E\in\pi(\mathcal B)$.
Let $X$ be an element in $\mathcal B$ such that
\begin{equation*}
	0\leqslant X\leqslant I,\qquad E=\pi(X).
\end{equation*}
Then $XG=GX=E$.
Thus, $XA=A$ and $XB=0$.
The cases $A=0$ or $B=0$ are immediate.
Therefore, $\mathcal A$ is an $SAW^*$-algebra.
\end{proof}

We conclude with several consequences of \Cref{thm sharply-clean}.
Recall that a $C^*$-algebra is \emph{$\sigma$-finite} if every family $\{A_j\}_{j\in\Lambda}$ of nonzero positive elements satisfying $A_iA_j=0$ for $i\ne j$ is countable; see \cite[Definition~2.1(b)]{HR-93}.
For a $C^*$-algebra $\mathcal B$, let $M(\mathcal B)$ denote its multiplier algebra and $M(\mathcal B)/\mathcal B$ its corona algebra.

\begin{corollary}\label{cor sharply-clean-consequences}
The following statements hold:
\begin{enumerate}
\item[$(1)$] A $\sigma$-finite unital $C^*$-algebra is sharply clean if and only if it is an $AW^*$-algebra.

\item[$(2)$] A separable unital $C^*$-algebra is sharply clean if and only if it is finite dimensional.

\item[$(3)$] Let $\mathcal B$ be a $\sigma$-unital $C^*$-algebra. The corona algebra $M(\mathcal B)/\mathcal B$ is sharply clean if and only if it has real rank zero.

\item[$(4)$] If $\{\mathcal A_n\}_{n=1}^\infty$ are unital $C^*$-algebras of real rank zero, then
\begin{equation*}
    \prod_{n=1}^\infty\mathcal A_n
    \bigg/
    \bigoplus_{n=1}^\infty\mathcal A_n
\end{equation*}
is sharply clean.
\item[$(5)$] If $\{\mathcal A_n\}_{n=1}^\infty$ are unital $C^*$-algebras of real rank zero and $\omega$ is a free ultrafilter on $\mathbb N$, then the norm ultraproduct
\begin{equation*}
    \prod_\omega\mathcal A_n
    =
    \prod_{n=1}^\infty\mathcal A_n
    \bigg/
    \left\{
        (A_n)_n\in\prod_{n=1}^\infty\mathcal A_n
        \colon \lim_{n\to\omega}\|A_n\|=0
    \right\}
\end{equation*}
is sharply clean.
\end{enumerate}
Here $\prod_{n=1}^\infty\mathcal A_n$ denotes the $C^*$-algebra of bounded sequences, and $\bigoplus_{n=1}^\infty\mathcal A_n$ is its ideal of sequences converging to zero in norm.
\end{corollary}

\begin{proof}
By \Cref{thm sharply-clean}, every sharply clean $C^*$-algebra is an $SAW^*$-algebra.
The necessity in $(1)$ follows from \cite[Section~3.1, p.~163]{HR-93}, and the necessity in $(2)$ follows from \cite[Corollary~2]{Pedersen-86}.
The converse implications follow from \Cref{prop AW}, since every finite-dimensional $C^*$-algebra is an $AW^*$-algebra.
This proves $(1)$ and $(2)$.

$(3)$ By \cite[Theorem~13]{Pedersen-86}, the corona algebra
$M(\mathcal B)/\mathcal B$ is an $SAW^*$-algebra.
The assertion therefore follows from \Cref{thm sharply-clean}.

$(4)$ Let $\mathcal B=\bigoplus_{n=1}^\infty\mathcal A_n$.
Then $\mathcal B$ is $\sigma$-unital and $M(\mathcal B)=\prod_{n=1}^\infty\mathcal A_n$.
Since the product has real rank zero and real rank zero passes to quotients, the conclusion follows from $(3)$.

$(5)$ Let $\mathcal I_\omega$ denote the ideal in the denominator in $(5)$.
Since $\omega$ is free, we have $\bigoplus_{n=1}^\infty\mathcal A_n\subseteq\mathcal I_\omega$.
Thus, $\prod_\omega\mathcal A_n$ is a quotient of the sharply clean algebra in $(4)$.
Since quotient maps preserve idempotents and invertibility and do not increase inverse norms, sharp cleanness passes to quotients.
Hence $\prod_\omega\mathcal A_n$ is sharply clean.
\end{proof}


\end{document}